  % FONTS

  \newcount\fontset
  \fontset=1
  \def\dualfont#1#2#3{\font#1=\ifnum\fontset=1 #2\else#3\fi}

  \dualfont\bbfive{bbm5}{cmbx5}
  \dualfont\bbseven{bbm7}{cmbx7}
  \dualfont\bbten{bbm10}{cmbx10}
  \font \eightbf = cmbx8
  \font \eighti = cmmi8 \skewchar \eighti = '177
  \font \eightit = cmti8
  \font \eightrm = cmr8
  \font \eightsl = cmsl8
  \font \eightsy = cmsy8 \skewchar \eightsy = '60
  \font \eighttt = cmtt8 \hyphenchar\eighttt = -1
  \font \msbm = msbm10
  \font \sixbf = cmbx6
  \font \sixi = cmmi6 \skewchar \sixi = '177
  \font \sixrm = cmr6
  \font \sixsy = cmsy6 \skewchar \sixsy = '60
  \font \tensc = cmcsc10
  
  \font \titlefont = cmr7 scaled \magstep4
  \scriptfont \bffam = \bbseven
  \scriptscriptfont \bffam = \bbfive
  \textfont \bffam = \bbten

  \font\rs=rsfs10 % For Lin(H)

  \newskip \ttglue

  \def \eightpoint {\def \rm {\fam0 \eightrm }%
  \textfont0 = \eightrm
  \scriptfont0 = \sixrm \scriptscriptfont0 = \fiverm
  \textfont1 = \eighti
  \scriptfont1 = \sixi \scriptscriptfont1 = \fivei
  \textfont2 = \eightsy
  \scriptfont2 = \sixsy \scriptscriptfont2 = \fivesy
  \textfont3 = \tenex
  \scriptfont3 = \tenex \scriptscriptfont3 = \tenex
  \def \it {\fam \itfam \eightit }%
  \textfont \itfam = \eightit
  \def \sl {\fam \slfam \eightsl }%
  \textfont \slfam = \eightsl
  \def \bf {\fam \bffam \eightbf }%
  \textfont \bffam = \eightbf
  \scriptfont \bffam = \sixbf
  \scriptscriptfont \bffam = \fivebf
  \def \tt {\fam \ttfam \eighttt }%
  \textfont \ttfam = \eighttt
  \tt \ttglue = .5em plus.25em minus.15em
  \normalbaselineskip = 9pt
  \def \MF {{\manual opqr}\-{\manual stuq}}%
  \let \sc = \sixrm
  \let \big = \eightbig
  \setbox \strutbox = \hbox {\vrule height7pt depth2pt width0pt}%
  \normalbaselines \rm }

  % HEADER

  \def \Headlines #1#2{\nopagenumbers
    \voffset = 2\baselineskip
    \advance \vsize by -\voffset
    \headline {\ifnum \pageno = 1 \hfil
    \else \ifodd \pageno \tensc \hfil \lcase {#1} \hfil \folio
    \else \tensc \folio \hfil \lcase {#2} \hfil
    \fi \fi }}

  \def \Title #1{\vbox{\baselineskip 20pt \titlefont \noindent #1}}

  \def \Date #1 {\footnote {}{\eightit Date: #1.}}

  \def \Authors #1{\bigskip \bigskip \noindent #1}

  \long \def \Addresses #1{\begingroup \eightpoint \parindent0pt
\medskip #1\par \par \endgroup }

  \long \def \Abstract #1{\begingroup \eightpoint
  \bigskip \bigskip \noindent
  {\sc ABSTRACT.} #1\par \par \endgroup }

  % CONTROL SEQUENCES

  \def \lcase #1{\edef \auxvar {\lowercase {#1}}\auxvar }
  \def \vg #1{\ifx #1\null \null \else
    \ifx #1\ { }\else
    \ifx #1,,\else
    \ifx #1..\else
    \ifx #1;;\else
    \ifx #1::\else
    \ifx #1''\else
    \ifx #1--\else
    \ifx #1))\else
    { }#1\fi \fi \fi \fi \fi \fi \fi \fi \fi }

  \def \goodbreak {\vskip0pt plus.1\vsize \penalty -250 \vskip0pt
plus-.1\vsize }

  \newcount \secno \secno = 0
  \newcount \stno

  \def \seqnumbering {\global \advance \stno by 1
    \number \secno .\number \stno }

  \def \label #1{\def\localvariable {\number \secno
    \ifnum \number \stno = 0\else .\number \stno \fi }\global \edef
    #1{\localvariable }}

  \def\section #1{\global\def\SectionName{#1}\stno = 0 \global
\advance \secno by 1 \bigskip \bigskip \goodbreak \noindent {\bf
\number \secno .\enspace #1.}\medskip \noindent \ignorespaces}

  \long \def \sysstate #1#2#3{\medbreak \noindent {\bf \seqnumbering
.\enspace #1.\enspace }{#2#3\vskip 0pt}\medbreak }
  \def \state #1 #2\par {\sysstate {#1}{\sl }{#2}}
  \def \definition #1\par {\sysstate {Definition}{\rm }{#1}}

  % Examples
  % \sysstate ...{Theorem}.{font}{Text}
  % \state .......Theorem.........Text\par
  % \definition ..................Text\par

  \def \proof {\medbreak \noindent {\it Proof.\enspace }}
  \def \proofend {\ifmmode \eqno \square \else \hfill \square
\looseness = -1 \medbreak \fi }

  \def \$#1{#1 $$$$ #1}
  \def\=#1{\buildrel (#1) \over =}

  \def\Item #1{\smallskip \item {#1}}
  \newcount \zitemno \zitemno = 0
  \def\izitem {\zitemno = 0}
  \def\zitem {\global \advance \zitemno by 1 \Item {{\rm(\romannumeral
\zitemno)}}}

  \newcount \footno \footno = 1
  \newcount \halffootno \footno = 1
  \def\footcntr {\global \advance \footno by 1
  \halffootno =\footno
  \divide \halffootno by 2
  $^{\number\halffootno}$}
  \def\fn#1{\footnote{\footcntr}{\eightpoint#1}}

  % STANDARD DEFINITIONS

  % STANDARD DEFINITIONS

  \def \({\left (}
  \def \){\right )}
  \def \[{\left \Vert }
  \def \]{\right \Vert }
  \def \*{\otimes }
  \def \+{\oplus }
  \def \:{\colon }
  \def \<{\left \langle }
  \def \>{\right \rangle }
  \def \text #1{\hbox {\rm #1}}
  \def \curly#1{\hbox{\rs #1\/}}
  \def \ds{\displaystyle}
  \def \and {\hbox {,\quad and \quad }}
  
  \def \calcat #1{\,{\vrule height8pt depth4pt}_{\,#1}}

  \def \crossproduct {{\hbox {\msbm o}}}
  
  \def \for #1{,\quad \forall\,#1}
  \def \inv {^{-1}}
  
  \def \square {\hbox {$\sqcap \!\!\!\!\sqcup $}}
  \def \stress #1{{\it #1}\/}

  \def \|{\Vert }
  \def \inv {^{-1}}

  % REFERENCE CONTROL

  \newcount \bibno \bibno =0
  \def \newbib #1{\global \advance \bibno by 1 \edef #1{\number
    \bibno}}
  \def\cite #1{{\rm [\bf #1\rm ]}}
  \def\scite #1#2{\cite{#1{\rm \hskip 0.7pt:\hskip 2pt #2}}}
  \def\lcite #1{(#1)}
  \def\fcite #1#2{\lcite{#1}}
  \def\bibitem#1#2#3#4{\smallskip \item {[#1]} #2, ``#3'', #4.}

  \def \references {
    \begingroup
    \bigskip \bigskip \goodbreak
    \eightpoint
    \centerline {\tensc References}
    \nobreak \medskip \frenchspacing }

  \def\AV{AV}
  \def\GAV{\widetilde{AV}}
  \def\CA{C^*}
  \def\B{\curly{B}}
  \def\G{\Gamma}
  \def\[{\big(}
  \def\]{\big)}
  \def\lS{\check S}
  \def\a{\alpha}
  \def\TCP{{\curly T}(A,\a,\Tr)}
  \def\N{{\bf N}}
  \def\Z{{\bf Z}}
  \def\T{{\bf T}}
  \def\U{{\curly U}}
  \def\K{{\cal K}}
  \def\E{{\cal E}}
  \def\M{{\cal M}}
  \def\I{\Lambda}
  \def\supp{{\rm supp}\,}
  
  \def\ind{{\rm ind}}
  \def\L{{\curly L}}
  \def\Tr{{\cal L}}
  \def\tS{{\hat S}}
  \def\d{\, d}
  
  \def\crpr{\mathop{\crossproduct_{\a,\Tr}} {\bf N}}
  \def\CP{A\crpr}
  \def\CCP{C(X)\crpr}
  \def\Li{L^\infty(X,\mu)}

  \def\fcite #1#2{#1}
  \def\yes{1}
  
  \def\showchanges{1}
  \def\newini{\medskip \line{\hrulefill\ BEGIN NEW \hrulefill} \medskip}
  \def\newini{\ifx\showchanges\yes
  \medskip \line{\hrulefill\ BEGIN NEW \hrulefill} \medskip
  \else\fi}
  \def\newend{\ifx\showchanges\yes
  \medskip \hrule \medskip
  \else \fi}

  \def \newbib #1#2{\global \advance \bibno by 1 \edef #1{\number
    \bibno}}

  \newbib{\alnr}{alnr}
  \newbib{\AS}{AS}
  \newbib{\ArzVershik}{ArzVershik}
  \newbib{\ArzVershikTwo}{ArzVershikTwo}
  \newbib{\Bourbaki}{Bourbaki}
  \newbib{\BKR}{BKR}
  \newbib{\Bratteli}{Bratteli}
  \newbib{\Deaconu}{Deaconu}
  \newbib{\EffrosHahn}{EffrosHahn}
  \newbib{\newpim}{newpim}
  \newbib{\Endo}{Endo}
  \newbib{\Tower}{Tower}
  \newbib{\RPF}{RPF}
  \newbib{\ELQ}{ELQ}
  \newbib{\anHuef}{anHuef}
  \newbib{\quasilat}{quasilat}
  \newbib{\Murphy}{Murphy}
  \newbib{\VN}{VN}
  \newbib{\Ped}{Ped}
  \newbib{\Power}{Power}
  \newbib{\Renault}{Renault}
  \newbib{\NewRenault}{NewRenault}
  \newbib{\Stacey}{Stacey}
  \newbib{\Takesaki}{Takesaki}
  \newbib{\Tom}{Tom}
  \newbib{\Watatani}{Watatani}

  \def\titletext{C*-Algebras of Irreversible Dynamical Systems}

  \Headlines {\titletext} {R.~Exel and A.~Vershik}

  \Title{\titletext}

  \Date{18 March 2002}

  \Authors
  {R.~Exel\footnote
    {*}{\eightrm Partially supported by CNPq.}
  and A.~Vershik\footnote
    {**}{\eightrm Partially supported by CRDF grant RM1-2244.}}

  \Addresses
  {Departamento de Matem\'atica,
  Universidade Federal de Santa Catarina,
  Florian\'opolis,
  Brazil
  (exel@mtm.ufsc.br).
  \par
  Russian Academy of Sciences, St.~Petersburg, Russia
  (vershik@pdmi.ras.ru).}

  \Abstract {We show that certain C*-algebras which have been studied
among others by Arzumanian, Vershik, Deaconu, and Renault in
connection to a measure preserving transformation of a measure space
and/or to a covering map of a compact space are special cases of the
endomorphism crossed-product construction recently introduced by the
first named author.  As a consequence these algebras are given
presentations in terms of generators and relations.  These results
come as a consequence of a general Theorem on faithfulness of
representations which are covariant with respect to certain circle
actions.  For the case of topologically free covering maps we prove a
stronger result on faithfulness of representations which needs no
covariance. We also give a necessary and sufficient condition for
simplicity.}

  \section{Introduction}
  Virtually all of the rich interplay between the theory of operator
algebras and Dynamical Systems occur by means of a construction
usually referred to as the \stress{crossed-product}.  The basic idea,
which was first employed by von Neumann \cite{\VN} in 1936 in the
context of automorphisms of measure spaces, consists in attaching an
operator  algebra (a W*-algebra in the case of von
Neumann's original construction) to a given dynamical
system whose algebraic structure is expected to reflect dynamical
properties of the given system.  The analog of this construction for
the case of C*-algebras is due to I.~Gelfand with co-authors (M.~Naimark,
S.~Fomin) who used it later for some special dynamical systems.
  
  By far the greatest advances in this enterprise have been achieved
for \stress{reversible} systems, i.e.~when the dynamics is implemented
by a group of invertible transformations.  This theory (sometimes
called algebraic theory of dynamical systems) now has many important
results and has become widespread and popular.
  Nevertheless a lot of
effort has been put into extending these advances to
\stress{irreversible} systems or \stress{dynamics of semigorups} and
thus an important theory emerged attempting to parallel the theory for
reversible systems.

The theory of C*-algebras for irreversible systems breaks down into
two somewhat disjoint areas, the first one consisting of the
construction and study of C*-algebras associated to classical
dynamical systems, i.e.~transformations of measure or topological
spaces (see
  e.g.~\cite{\ArzVershik},
  \cite{\ArzVershikTwo},
  \cite{\Deaconu},
  \cite{\NewRenault}), while the second deals with the study of
crossed-products of C*-algebras by endomorphisms or semigroups thereof
(see
  e.g.~\cite{\Stacey},
  \cite{\BKR}, 
  \cite{\alnr},
  \cite{\quasilat}, 
  \cite{\Murphy}).

  As it turned out these two areas have had little intersection partly
because the latter, which perhaps received the greatest share of
attention in recent times, boasts its greatest successes when the
endomorphisms involved have a hereditary range, a hypothesis which is
absent when one deals with endomorphisms arising from classical
systems.

  The difference between the C* and W*-algebra approach is more
serious in the case of irreversible dynamics than in the reversible
one.  The point is that the involution in the C*-algebra for
irreversible case must be defined separately meanwhile for the group
case the involution is defined more or less uniquely. At the same time
if we consider measure preserving dynamics and $L^2$-theory the
involution exists automatically and this gives possibility to define
also C*-algebras.  This was the way of the papers \cite{\ArzVershik},
\cite{\ArzVershikTwo} --- to consider measure preserving action of the
endomorphism and the C*-hull in $L^2$ of the multiplicators and
endomorphism. After this we can define another representation of
this C*-algebra which leads to von Neumann factors.  Actually we do
not need the invariant measure itself but only the conditional expectation
(e.g. conditional measure on the partition onto pre-images of the
points).  In pure topological situation we can firstly consider very
large C*-algebra and then factorize it using some conditional
expectation (see below).

In this paper our interest is on the oldest, and perhaps not so
intensively studied, relationship between classical irreversible
systems and C*-algebras.  Our focus is on two classes of C*-algebras,
the first one being defined in terms of a measure preserving
transformation $T$ of a measure space $(X,\mu)$, and is closely
related to the Arzumanian--Vershik algebras (see \cite{\ArzVershik}).
These algebras, which are defined below, will be denoted here by
$\GAV(X,T,\mu)$.

The second class of C*-algebras that we study
here arises from the consideration of self covering maps $T:X\to X$,
where $X$ is a compact space  (see \cite{\ArzVershikTwo}, \cite{\Deaconu}, and
\cite{\NewRenault}).  The algebras in this second class will be 
denoted by $\CA(X,T)$.

It is the main goal of this work to give a unified treatment of these 
algebras by first showing 
  (Theorems \fcite{6.1}{\ArzVershikIsomCrossProd} and 
  \fcite{9.1}{\CoveringIsomCrossProd} below)
that they are special cases of the
endomorphism crossed-product construction recently introduced by the
first named author (incidentally that construction also includes crossed-products
by endomorphisms with hereditary range \scite{\Endo}{Section 4}). 

In addition we show how to use this unified perspective to prove new
results about $\GAV(X,T,\mu)$ and $\CA(X,T)$
  % as special cases of general results for endomorphism
  % crossed-products, some of which relate to faithful representations
  along the lines of
  \cite{\Stacey},
  \cite{\BKR},
  \cite{\alnr},
  \cite{\quasilat},
  \cite{\ELQ}.
  First we prove a version of what has been called the ``gauge
invariant uniqueness Theorem'' \scite{\anHuef}{2.3} by showing
  (Theorem \fcite{4.2}{\GaugeInvariantUniqueness})   
  that any representation of the crossed-product which is faithful on
the core algebra and which is covariant relative to certain circle
actions must be faithful on all of the crossed-product algebra.  The
application of this to either $\AV(X,T,\mu)$ or $\CA(X,T)$ lead to
Theorems on faithfulness of representations which are perhaps hitherto
unknown.

As another application we show
  (Theorem \fcite{5.5}{\FaithfulRepFromPhi})
  how the proposal outlined in the closing paragraph of
\cite{\ArzVershik} to extend the definition of the Arzumanian--Vershik
algebra to a non-commutative setting can be modified to achieve a down
to the ground elementary construction of the crossed-product algebra
(assuming the existence of certain faithful states).

We then specialize to the case in which the range of the endomorphism
considered admits a conditional expectation of \stress{index-finite
type} according to Watatani \cite{\Watatani}.  We first show
  (Corollary \fcite{7.3}{\ShorterRelations})
  that the process of ``modding out the redundancies'' in Definition
(\fcite{2.7}{\DefCrossProd}) below may be achieved by introducing a
single relation and hence we obtain a straightforward presentation of
the crossed-product algebra in terms of generators and relations.

The situation leading up to $\CA(X,T)$, namely that of a covering map,
is shown to fall under the hypothesis of finiteness of the index and
hence we obtain
  (Theorem \fcite{9.2}{\ShorterRelationsForCovering})
  a very concise set of relations defining $\CA(X,T)$.

Still in the case of a covering map we obtain another result on
faithfulness of representations,  this time without regard for circle
actions but under the hypothesis of \stress{topological
freeness}: a continuous map $T:X\to X$ is said to be topologically
free\fn{Also called essentially free in \cite{\Deaconu}.} (see also
  \scite{\Tom}{2.1},
  \cite{\AS},
  \cite{\Deaconu}, and
  \scite{\ELQ}{2.1})
  if for every pair of nonnegative integers $(n,m)$ with $n\neq m$ one
has that the set
  $
  \{ x \in X : T^n(x) = T^m(x) \}
  $
  has empty interior.
  Precisely we show
  (Theorem \fcite{10.3}{\MainResultTopFree})
  that under this hypothesis any representation of
the crossed-product which is faithful on $C(X)$ must itself be faithful.
  We moreover show
  (Theorem \fcite{11.2}{\SimplicityResult})
  that the crossed-product algebra is simple if and
only if $T$ is irreducible
  in the sense that there are no closed nontrivial sets $F\subseteq X$
such that $T\inv(F)=F$.  This result improves on Deaconu's
characterization of simplicity \cite{\Deaconu} since the above notion
of irreducibility is much weaker than Deaconu's notion of minimality.

The hypothesis that $T$ is a covering map is directly related to the
finiteness of the index which in turn is responsible for the existence
of a conditional expectation from the crossed-product to $C(X)$
\scite{\Tower}{8.9}.  Since that conditional expectation turns out to
be one of our main tools in proving the characterization of simplicity
our methods completely break down in the infinite index case.  The
characterization of simplicity in the infinite index case therefore
stands out as one of the main questions left unresolved.

Last but not least we should mention that among the algebras which can
be described as $\CCP$ are the Cuntz--Krieger algebras
\scite{\Endo}{6.2} and the Bunce--Deddens algebras \cite{\Deaconu},
which therefore provide examples of the algebras that we discuss here.

Part of this work was done while the authors  were visiting
  the Centre for Advanced Study
  at the Norwegian Academy of Science and Letters.
  We would like to thank Ola Bratteli and 
  Magnus Landstad, 
  the organizers there, for their warm hospitality.

  \section{The algebras}
  Perhaps the first attempt at constructing an operator algebra from an
irreversible dynamical system is Arzumanian and Vershik's 1978 paper
\cite{\ArzVershik}.  Given a (not necessarily invertible) measure
preserving transformation $T$ of a measure space $(X,\mu)$
  % Arzumanian and Vershik
  they
  consider the C*-algebra of operators on $L^2(X,\mu)$ generated by
the \stress{multiplication operators}
  $$
  M_f : \xi\in L^2(X,\mu) \ \longmapsto \ f\xi \in L^2(X,\mu),
  $$
  for all $f\in \Li$, in addition to the isometry $S$ defined by
  $$
  S(\xi)\calcat x = \xi\big(T(x)\big)
  \for \xi\in L^2(X,\mu)
  \for x\in X.
  $$
  We will refer to this algebra as the \stress{Arzumanian--Vershik}
algebra and we will denote it by $\AV(X,T,\mu)$.

  While it would be desirable that the algebraic structure of
$\AV(X,T,\mu)$ depend more on the dynamical properties of $T$ than on
the invariant measure $\mu$ chosen it is evident that it does depend
heavily on $\mu$.  In fact there is a specific way in which $\mu$
influences the algebraic structure of $\AV(X,T,\mu)$ which we would
like to describe.  To this goal assume that we are under the
hypothesis of \scite{\Bourbaki}{\S3, N$^o$1, Th\'eor\`eme 1} so that
$\mu$ may be disintegrated along the fibers of $T$.  This means that
there exists a collection $\{\mu^x\}_{x\in X}$, where each $\mu^x$ is a
probability measure on $T\inv(x)$, and for every
$\mu$-integrable function $f$ on $X$ one has that
  $$
  \int_X f(x)\d \mu(x) =
  \int_X \(\int_{T\inv(x)} f(y)\d \mu^x(y)\) \d \nu(x),
  \eqno {(\seqnumbering)}
  \label \Disintegration
  $$
  where $\nu = T^*(\mu)$ (observe that $\nu = \mu$ since $\mu$
is invariant under $T$).  For every $f\in \Li$ define the function $\Tr(f)$ on
$X$ by
  $$
  \Tr(f)\calcat x = \int_{T\inv(x)} f(y) \d \mu^x(y).
  \eqno {(\seqnumbering)}
  \label \VershikTransfer
  $$
  One may then verify that
  $$
  S ^* M_f S = M_{\Tr(f)},
  \eqno {(\seqnumbering)}
  \label \VershikRelation
  $$
  for every $f\in \Li$.  One therefore sees that the measure $\mu$, or
at least the collection formed by the $\mu^x$, does indeed influence
the algebraic structure of $\AV(X,T,\mu)$.

Observe that in case $T$ is the identity map on $X$ one has that $S$
is the identity operator on $L^2(X,\mu)$.  Because this is not in
accordance with the classical case of crossed-products by
automorphisms (where the unitary element implementing the given
automorphism has full spectrum) it is sometimes useful to change the
definition of $\AV(X,T,\mu)$ as follows:  letting $U$ be the bilateral
shift on $\ell^2(\Z)$ define the operators $\tilde S$ and $\tilde M_f$
on the Hilbert space $L^2(X,\mu)\* \ell^2(\Z)$ by
  $$
  \tilde S = S\*U
  \and
  \tilde M_f = M_f\*1,
  $$
  for all $f\in \Li$. We then let $\GAV(X,T,\mu)$ be the C*-algebra of
operators on $L^2(X,\mu)\* \ell^2(\Z)$ generated by $\tilde S$ and all
the $\tilde M_f$.
  The advantage of this alternate construction is that when $T$ is an
automorphism of $X$ one has that
  $$
  \GAV(X,T,\mu) \simeq \Li \crossproduct_T \Z,
  $$
  where the right hand side is the classical crossed-product
\cite{\Ped} of $\Li$ under the automorphism induced by $T$.

  % It is clear that the obvious generalization of relation
  % \lcite{\VershikRelation} holds in $\GAV(X,T,\mu)$.

There is another important construction which has been proposed as a
replacement to classical crossed-products in case $X$ is a compact
space and $T:X\to X$ is a covering map.  Building on
\cite{\ArzVershikTwo} Deaconu \cite{\Deaconu} has shown how to
construct an $r$-discrete groupoid from the pair $(X,T)$ to which one
may apply Renault's well known construction \cite{\Renault} leading up
to a C*-algebra which we will denote by $\CA(X,T)$.

  Still under the assumption that $T$ is a covering map 
  % (irrespective of $\#T\inv(x)$ being constant or not)
  each inverse image $T\inv(x)$ is a
finite set which in turn admits a canonical probability measure
$\mu^x$, namely the normalized counting measure.  Expression
\lcite{\VershikTransfer} may then be employed and this time it defines
a bounded positive operator
  $$
  \Tr: C(X) \to C(X).
  $$
  Considering $C(X)$ as a subalgebra of $\CA(X,T)$ and letting $S$ be
the isometry in $\CA(X,T)$ described shortly before the statement of
Theorem 2 in \cite{\Deaconu} it is easy to see that
  $$
  S ^* f S = \Tr(f)
  \for f\in C(X),
  $$
  which is seen to be another manifestation of
\lcite{\VershikRelation}.

  In a recent article \cite{\Endo} the first named author has proposed
yet another construction of crossed-product in the irreversible
setting (see also \cite{\Tower}, \cite{\RPF}), the ingredients of
which are a unital C*-algebra $A$, an injective *-endomorphism
  $$
  \a : A \to A
  $$
  such that $\a(1)=1$, and a \stress{transfer operator}, namely a
continuous positive linear map
  $$
  \Tr: A \to A
  $$
  such that
  $$
  \Tr\(a\a(b)\) = \Tr(a)b
  \for a,b\in A,
  \eqno {(\seqnumbering)}
  \label \DefTransfer
  $$
  which we will moreover suppose satisfies $\Tr(1)=1$.

The main examples we have in mind are of course related to the above
situations, where $A$ is either $\Li$ or $C(X)$ and $\a$ is given in
both cases by the expression
  $$
  \a: f \in A \mapsto f\circ T \in A.
  \eqno {(\seqnumbering)}
  \label \CommutativeEndo
  $$
  As for a transfer operator one may use \lcite{\VershikTransfer}
after the appropriate choice of probability measures $\mu^x$ is made.

Let us fix for the time being a C*-algebra $A$, an injective
*-endomorphism $\a$ of $A$ preserving the unit, and a transfer
operator $\Tr$ such that $\Tr(1)=1$.

  \definition
  \label \DefTopExt
  \scite{\Endo}{3.1}\quad
  We will let $\TCP$ be the universal unital C*-algebra generated by a
copy of $A$ and an element $\tS$ subject to the relations:
  \izitem
  \zitem $\tS a = \a(a)\tS $, and
  \zitem $\tS^*a\tS = \Tr(a)$,
  \medskip\noindent
  for every $a\in A$.

  As proved in \scite{\Endo}{3.5} the canonical map from $A$ to $\TCP$
is injective so we may view $A$ as a subalgebra of $\TCP$.  Since
$\Tr(1)=1$ we have that
  $$
  \tS^*\tS = \tS^*1\tS = \Tr(1) = 1,
  $$
  and hence we see that $\tS $ is an isometry.

Following \scite{\Endo}{3.6} a \stress{redundancy} is a pair $(a,k)$
of elements in $\TCP$ such that $k$ is in the closure of $A\tS
\tS^*A$, $a$ is in $A$, and
  $$
  ab\tS = kb\tS \for b\in A.
  $$

  \definition
  \label \DefCrossProd
  \scite{\Endo}{3.7} \quad
  The crossed-product of $A$ by $\a$ relative to $\Tr$, denoted by
$\CP$, is defined to be the quotient of $\TCP$ by the closed two-sided
ideal generated by the set of differences $a-k$, for
  all\fn{We should remark that in Definition 3.7 of \cite{\Endo} one
uses only the redundancies $(a,k)$ such that $a$ belongs to the ideal
of $A$ generated by the range of $\a$.  But, under the present
hypothesis that $\a$ preserves the unit, this ideal coincides with
$A$.}
  redundancies $(a,k)$.

  As stated in \scite{\Tower}{4.12} the canonical map from $A$ to
$\CP$ is injective and we use it to think of $A$ as a subalgebra of
$\CP$.
  We will denote by
  % $q$ the canonical quotient map
  % $$
  % q: \TCP \to \CP,
  % % \eqno {(\seqnumbering)}
  % % \label \NotationQuotient
  % $$
  % and by
  $S$ the image of $\tS$ in $\CP$.

  \section{Preliminaries}
  For  ease of reference let us now list the assumptions
which will be in force for the largest part of this work.

  \sysstate{Standing Hypotheses}{\rm}{\label \StandingHyp
  \izitem
  \zitem $A$ is a unital C*-algebra,
  \zitem $\a$ is an injective *-endomorphism of $A$ such that
$\a(1)=1$,
  \zitem $\Tr$ is a normalized transfer operator, i.e.~a continuous
positive linear map from $A$ to itself such that $\Tr(1)=1$, and
  $
  \Tr\(a\a(b)\) = \Tr(a)b,
  $
  for all $a$ and $b$ in $A$ 
  (see \scite{\Endo}{Definition 2.1}).}

  We now wish to review some results from \cite{\Tower} in a
form which is suitable for our purposes.  Recall from
\scite{\Tower}{2.3} that $\TCP$ (resp.~$\CP$) coincides with the
closed linear span of the set of elements of the form
$a\tS^n\tS^{*m}b$ (resp.~$aS^nS^{*m}b$), for $a,b\in A$.

  \state Proposition
  \label \ConditionalExpectation
  {\rm (See also \scite{\Endo}{2.3}).}
  For each $n\in\N$ the map $\E_n$ defined by
  $$
  \E_n=\a^n\circ\Tr^n
  $$
  is a conditional expectation from $A$ onto the range of $\a^n$.

  \proof
  It is clear that the range of $\E_n$ is contained in the range of
$\a^n$.  Moreover if $a\in A$ and $b\in\a^n(A)$, say $b = \a^n(c)$
with $c\in A$, we have that
  $$
  \E_n(ab) =
  \a^n\big(\Tr^n(a\a^n(c))\big) =
  \a^n\big(\Tr^n(a)\big) \a^n(c) =
  \E_n(a)b.
  $$
  Plugging $a=1$ above, and noticing that $\E_n(1)=1$, we see that
$\E_n$ is the identity on the range of $\a^n$.  The remaining details
are now elementary.
  \proofend

  \definition
  We will denote by $\M_n$ and $\K_n$ the subsets of $\CP$ given by
  $$
  \M_n=\overline{AS^n}
  \and
  \K_n = \overline{\M_n\M_n^*} = \overline{AS^nS^{*n}A}.
  $$

  \state Proposition
  \label \MisModuleOverK
  For every $n,m\in\N$ with $n\leq m$ we have that
  \izitem
  \zitem $\K_n \M_m \subseteq \M_m$, and
  \zitem $\K_n \K_m \subseteq \K_m$.

  \proof Given $a,b,c\in A$ we have
  $$
  (a S^nS^{n*}b)(c S^m) =
  a S^nS^{n*} b c S^n S^{m-n} =
  a \a^n(\Tr^n(bc)) S^m,
  $$
  from where (i) follows.  As for (ii) we have
  $$
  \K_n \K_m \subseteq
  \overline{\K_n \M_m\M_m^*} \subseteq
  \overline{\M_m\M_m^*} = \K_m.
  \proofend
  $$

  It follows that each $\K_n$ is a C*-subalgebra of $\CP$.  In some
cases (see below) the $\K_n$ form an increasing sequence.
In any event the closure of the sum of all $\K_n$'s is a
C*-subalgebra which we denote by
  $$
  \U =
   \overline{\vrule height 10pt width 0pt\hbox{$\sum_{n=0}^\infty$}
\K_n}.
  $$
  Obviously $\U$ is the linear span of the set of elements of the form
  $aS^nS^{*n}b$, where $n\in \N$ and $a,b\in A$.

  There is
  (see \scite{\Tower}{3.3}) a canonical action $\gamma$ of the circle
group $\T$ on $\CP$ given by
  $$
  \gamma_z(S) = zS
  \and
  \gamma_z(a) = a
  \for a\in A
  \for z\in \T,
  $$
  which we shall call the \stress{gauge action}.
  For any circle action the expression
  $$
  F(a) = \int_{z\in \T} \gamma_z(a) \d z
  \for a\in\CP
  \eqno {(\seqnumbering)}
  \label \CondExpF
  $$
  gives a faithful conditional expectation onto the algebra of fixed
points for $\gamma$.  By \scite{\Tower}{3.5} we have that this fixed
point algebra is precisely $\U$.

  \section {Faithfulness of covariant representations}
  \label \SectionFaithfulness
  Our major goal in this section will be to establish sufficient
conditions for a representation of $\CP$ to be faithful.  We therefore
fix throughout this section a *-homomorphism
  $$
  \pi: \CP \to B,
  $$
  where $B$ is any C*-algebra, quite often chosen to be the algebra of
all bounded operators on a Hilbert space, in which case we will be
talking of a representation of $\CP$ as indicated above.

  \state Lemma
  \label \faithfulOnU
  If $\pi$ is faithful on $A$ then $\pi$ is faithful on all of $\U$.

  \proof We begin by claiming that $\pi$ is isometric when restricted
to $\M_n$ for all $n\in \N$.  In order to see this let $a\in A$ and
notice that
  $$
  \|\pi(aS^n)\|^2 =
  \|\pi(S^{*n}a^*aS^n)\| =
  \|\pi(\Tr^n(a^*a))\| =
  \|\Tr^n(a^*a)\| =
  \|S^{*n}a^*aS^n\| =
  \|aS^n\|^2,
  $$
  hence proving the claim.

  Let $k \in \U$ and suppose that $k = \sum_{i=0}^n k_i$, where each
$k_i\in \K_i$.  By \lcite{\MisModuleOverK.i} we have that $kbS^n \in
\M_n$ for every $b\in A$.
  Supposing that $\pi(k)=0$ we have that $\pi(kbS^n)=0$ for every $b$
in $A$ and hence $kbS^n=0$ by the above claim.  It follows that
  $kk'=0$ for all $k'\in\K_n$ and hence that $(k_0,k_1,\ldots,k_n)$ is
a redundancy of order $n$ as defined in \scite{\Tower}{6.2}.  By
\scite{\Tower}{6.3} we have that $k=0$ in $\CP$ and hence that $\pi$
is injective on $\sum_{i=0}^n \K_i$.  Since the latter is a C*-algebra
by \lcite{\MisModuleOverK.ii} and \scite{\Ped}{1.5.8} we see that
$\pi$ is isometric there and hence also on the inductive limit of
these algebras as $n\to\infty$ which coincides with $\U$.
  \proofend

This allows us to prove a version of the ``gauge invariant uniqueness
Theorem'' as in \scite{\anHuef}{2.3}.

  \state Theorem
  \label \GaugeInvariantUniqueness
  Under the assumptions in \lcite{\StandingHyp} let $B$ be a
C*-algebra and $\pi: \CP \to B$ be a *-homomorphism which is faithful
on $A$.  Suppose that $B$ admits an action of the circle group such
that $\pi$ is covariant relative to the gauge action on $\CP$ and the
given action on $B$.  Then $\pi$ is faithful on all of $\CP$.

  \proof
  By \scite{\newpim}{2.9} if suffices to prove that $\pi$ is faithful
on the fixed point algebra for the gauge action, which we know is
$\U$.  But this obviously follows from \lcite{\faithfulOnU}.
  \proofend

  %--------------
  \section{States}
  In this section we will develop the necessary tools in order to show
that the generalized Arzumanian--Vershik algebra $\GAV(X,T,\mu)$ fits
our crossed-product construction.  Since we will be working in the
non-commutative setting what we do here, especially Theorem
(\fcite{5.5}{\FaithfulRepFromPhi}) below, is related to the situation
discussed in the last paragraph of \cite{\ArzVershik}.  As before we
keep \lcite{\StandingHyp} in force.

  Throughout this section we will fix a state $\phi$ on $A$ which is
invariant under both $\a$ and $\Tr$, that is,
  $$
  \phi(\a(a)) = \phi(a) = \phi(\Tr(a))
  \for a\in A.
  $$
  Observe that $\phi$ is then necessarily invariant under the
conditional expectations $\E_n$ of \lcite{\ConditionalExpectation}.

  Let $\rho$ be the GNS representation associated to $\phi$ and denote
by $H$ the corresponding Hilbert space and by $\xi$ the cyclic vector.
Because $\phi$ is $\a$-invariant one has that the correspondence
  $$
  \rho(a)\xi \mapsto \rho\[\a(a)\big)\xi
  \for a\in A
  $$
  preserves norm and hence extends to $H$ giving an isometry $\lS\in
\B(H)$ such that
  $$
  \lS\[\rho(a)\xi\] = \rho\[\a(a)\]\xi
  \for a\in A.
  \eqno {(\seqnumbering)}
  \label \RepresentedS
  $$
  We next claim that for every $a\in A$ one has that
  $$
  \lS^*\[\rho(a)\xi\] = \rho\[\Tr(a)\]\xi.
  \eqno {(\seqnumbering)}
  \label \ClaimAdjoint
  $$
  In fact, given $b\in A$ we have
  $$
  \<\lS^*\[\rho(a)\xi\], \rho(b)\xi\> =
  \<\rho(a)\xi, \lS\[\rho(b)\xi\]\> =
  \<\rho(a)\xi, \rho\[\a(b)\]\xi\> =
  \phi\[\a(b^*)a\] \$=
  \phi\[\Tr(\a(b^*)a)\] =
  \phi\[b^*\Tr(a)\] =
  \<\rho\[\Tr(a)\]\xi, \rho(b)\xi\>,
  $$
  from which \lcite{\ClaimAdjoint} follows.
  With this it is easy to show that
  for every $a\in A$ one has that
  $$
  \lS\rho(a) = \rho\[\a(a)\]\lS
  \and
  \lS^*\rho(a)\lS = \rho\[\Tr(a)\].
  $$
  By the universal property it follows that there exists a
representation $\Pi$ of $\TCP$ on $H$ such that
  $$
  \Pi(\tS) = \lS
  \and
  \Pi(a) = \rho(a),
  $$
  for all $a\in A$.

  \state Proposition
  \label \RepFromPhi
  For any redundancy $(a,k)$ one has that $\Pi(a-k)=0$.  Therefore
$\Pi$ factors through a representation $\pi$ of $\CP$ on $H$ such that
$\pi(S) = \lS$ and $\pi(a) = \rho(a)$ for all $a\in A$.

  \proof
  Plugging $a=1$ in \lcite{\RepresentedS} we see that $\lS(\xi) = \xi$
and hence for any $b\in A$ we have
  $$
  \Pi(b\tS)\xi =
  \rho(b)\lS\xi = \rho(b)\xi.
  $$
  Therefore
  $$
  \rho(a)\rho(b)\xi=
  \rho(ab)\xi =
  \Pi(ab\tS)\xi =
  \Pi(kb\tS)\xi =
  \Pi(k)\Pi(b\tS)\xi =
  \Pi(k)\rho(b)\xi,
  $$
  which implies that $\rho(a) = \Pi(k)$.
  \proofend

There is another representation of $\CP$ obtained from $\phi$ which we
would like to describe next.

  \state Proposition
  \label \RepFromPhiAndU
  Let $\rho$ be the above GNS representation of $A$ and let $U$ be the
bilateral shift on $\ell^2(\Z)$.  Then there exists a representation
$\tilde\pi$ of $\CP$ on $H\*\ell^2(\Z)$ such that
  $$
  \tilde\pi(S) = \lS\*U
  \and
  \tilde\pi(a) = \rho(a)\*1,
  $$
  for all $a\in A$.

  \proof
  It is elementary to check that
  $$
  (\lS\*U)(\rho(a)\*1) = \big(\rho\[\a(a)\]\*1\big)(\lS\*U)
  $$ and $$
  (\lS\*U)^*(\rho(a)\*1)(\lS\*U) = \rho\[\Tr(a)\]\*1,
  $$
  from where one concludes that there exists a representation
$\tilde\Pi$ of $\TCP$ on $H\*\ell^2(\Z)$ such that
  $$
  \tilde\Pi(\tS) = \lS\*U
  \and
  \tilde\Pi(a) = \rho(a)\*1.
  $$
  As before we claim that $\tilde\Pi(a-k)=0$ for any redundancy
$(a,k)$.  In order to see this observe that given $b,c\in A$ we have
  $$
  \tilde\Pi(b \tS \tS^* c) =
  (\rho(b)\*1)
  (\lS\*U)
  (\lS^*\*U^*)
  (\rho(c)\*1) =
  \Pi(b \tS \tS^* c) \* 1.
  $$
  This implies that
  $\tilde\Pi(k) = \Pi(k) \* 1$ for all $k\in \overline{A \tS \tS^*A}$
and hence the claim follows from \lcite{\RepFromPhi}.  It therefore
follows that $\tilde\Pi$ factors through $\CP$ yielding the desired
representation.
  \proofend

  The main result of this section is in order:

  \state Theorem
  \label \FaithfulRepFromPhi
  Under \lcite{\StandingHyp} let $\phi$ be a state on $A$ which is
invariant under both $\a$ and $\Tr$.  Suppose that the GNS
representation of $\phi$, denoted $\rho$, is faithful.  Then the
representation $\tilde\pi$ of $\CP$ constructed above is faithful.
Therefore $\CP$ is isomorphic to the algebra of operators on
$H\*\ell^2(\Z)$ generated by the set
  $$
  \{\rho(a)\*1 : a\in A\} \cup \{\lS\*U\},
  $$
  where $\lS$ is given by \lcite{\RepresentedS} and $U$ is the
bilateral shift on $\ell^2(\Z)$.

  \proof
  For each $z\in \T$ let $V_z$ be the unitary operator on $\ell^2(\Z)$
given on the canonical basis $\{e_n\}_{n\in \Z}$ by $V_z(e_n) = z^n
e_n$.  We then have that $V_z U V_z^* = z U$ which implies that
  $$
  (1\*V_z) \tilde\pi(S) (1\*V_z^*) = z\ \tilde\pi(S),
  $$
  while
  $$
  (1\*V_z) \tilde\pi(a) (1\*V_z^*)= \tilde\pi(a),
  $$
  for all $a\in A$.
  Let $B$ be the C*-algebra formed by the operators $T$ on
$H\*\ell^2(\Z)$ for which the map $z \mapsto (1\*V_z) T (1\*V_z^*)$ is norm
continuous.  The above calculations show that $\tilde\pi$ maps $\CP$
into $B$.  Moreover it is clear that $\tilde\pi$ is covariant for the
dual action on $\CP$ and the action by conjugation by $1\*V_z$ on $B$.
  The conclusion then follows from \lcite{\GaugeInvariantUniqueness}.
  \proofend

  \section{The Arzumanian--Vershik algebra}
  Let $(X,\mu)$ be a probability space and $T:X\to X$ be a measure
preserving transformation.  It is easy to see that the correspondence
  $$
  \a: f \in \Li \mapsto f\circ T \in \Li
  $$
  gives an injective *-endomorphism of $\Li$.
  We will assume that $\mu$ may be disintegrated along the fibers of
$T$ as in \lcite{\Disintegration}.  Letting $\{\mu^x\}_{x\in X}$ be
such a disintegration it is straightforward to verify that the
operator $\Tr$ defined by \lcite{\VershikTransfer} is a normalized
transfer operator.
  We therefore have all of the ingredients needed to form the
crossed-product algebra $\Li\crpr$.

  Let $\phi$ be the state on $\Li$ given by integration against $\mu$.
It is clear that $\phi$ is invariant under $\a$.  On the other hand
observe that \lcite{\Disintegration} together with the fact mentioned
there that $\nu=\mu$ says precisely that $\phi$ is invariant under
$\Tr$ as well.

  It is easy to see that $\AV(X,T,\mu)$ is precisely the range of the
representation $\pi$ of \lcite{\RepFromPhi} while $\GAV(X,T,\mu)$ is
the range of the representation $\tilde\pi$ of
\lcite{\RepFromPhiAndU}.  Since it is obvious that the GNS
representation associated to $\phi$ is faithful on $\Li$ we conclude
from \lcite{\FaithfulRepFromPhi} that:

  \state Theorem
  \label \ArzVershikIsomCrossProd
  The generalized Arzumanian--Vershik algebra $\GAV(X,T,\mu)$ is
isomorphic to the crossed-product $\Li\crpr$.

The above result does not rule out the possibility that $\AV(X,T,\mu)$
be isomorphic to $\Li\crpr$ as well.  Since the former algebra is a
quotient of the latter by \lcite{\RepFromPhi}, these algebras would be
isomorphic in case $\Li\crpr$ is simple, for example.

  Unfortunately we do not have a precise characterization of
simplicity, but see below for a similar characterization when $T$ is a
covering map of a compact space $X$.

Theorem \lcite{\ArzVershikIsomCrossProd} also helps explain to which
extent does $\GAV(X,T,\mu)$ depend on $\mu$.  While there is no
question that the probability measures $\mu^x$ influence the outcome
we see that this is as far as this dependence goes.  In particular if
one takes another measure $\nu$ such that the $\nu^x$ coincide with
the $\mu^x$ the algebras $\GAV(X,T,\mu)$ and $\GAV(X,T,\nu)$ will be
isomorphic.

  \section{Finite index endomorphisms}
  In this section we return to the general case described by
\lcite{\StandingHyp}.
  We will denote the the conditional expectation $\E_1$ of
\lcite{\ConditionalExpectation} simply by $E$.
  The main assumption to be added from here on is that $E$ is of
\stress{index-finite type} according to \scite{\Watatani}{1.2.2 and
2.1.6}.  That is, we will assume the existence of a
\stress{quasi-basis} for $E$, namely a finite sequence
$\{u_1,\ldots,u_m\}\subseteq A$ such that
  $$
  a=\sum_{i=1}^m u_i E(u_i^* a)\for a\in A.
  \eqno {(\seqnumbering)}
  \label\DefFiniteIndex
  $$
  In this case one defines the \stress{index} of $E$ by
  $$
  \ind(E) = \sum_{i=1}^m u_iu_i^*.
  $$
  It is well known that $\ind(E)$
  does not depend on the choice of the $u_i$'s,
  that it belongs to the center of $A$ \scite{\Watatani}{1.2.8}
  and
  is invertible \scite{\Watatani}{2.3.1}.

  In particular $A$ is a finitely generated right $B$-module with
$\{u_1,\ldots,u_m\}$ being a generating set.  See also
(\fcite{8.6}{\CoveringIsFinite}) below for a sufficient condition for $E$ to
be of index-finite type in the case that $A$ is a commutative C*-algebra.

  We would first like to show that under the present hypothesis the
process of modding out the redundancies in the definition of $\CP$
\lcite{\DefCrossProd} can be achieved by adding a single relation.

  \state Proposition
  \label \MainRedundancy
  Let $\{u_1,\ldots,u_m\}$ be a quasi-basis for $E$ and denote by
$k_0$ the element of $\TCP$ given by
  $
  k_0 = \sum_{i=1}^m u_i\tS\tS^*u_i^*.
  $
  Then the pair $(1,k_0)$ is a redundancy and hence
  $$
  1= \sum_{i=1}^m u_iSS^*u_i^*
  $$
  in $\CP$.  Moreover the kernel of the natural quotient map $q : \TCP
\to \CP$ coincides with the ideal generated by $1 - k_0.$

  \proof Given $b\in A$ we have
  $$
  k_0b\tS =
  \sum_{i=1}^m u_i\tS\tS^*u_i^* b\tS =
  \sum_{i=1}^m u_i \a(\Tr(u_i^*b))\tS =
  \sum_{i=1}^m u_i E(u_i^*b)\tS =
  b \tS,
  $$
  proving the first statement.  Next let $(a,k)$ be any other
redundancy.  Then
  $$
  k k_0 =
  \sum_{i=1}^m k u_i\tS\tS^*u_i^* =
  \sum_{i=1}^m a u_i\tS\tS^*u_i^* =
  a k_0.
  $$
  It follows that
  % I like this formula displayed !
  $$
  a - k = (a-k)(1-k_0)
  $$
  % I like this formula displayed !
  and hence that $a-k$ belongs to the ideal generated by $1-k_0$.
  \proofend

  The description of $\CP$ in terms of generators and relations
therefore becomes:

  \state Corollary
  \label \ShorterRelations
  Under \lcite{\StandingHyp} suppose that the conditional expectation
$E=\a\circ\Tr$ is of index-finite type and let 
$\{u_1,\ldots,u_m\}$ be a quasi-basis for $E$.  Then
  $\CP$ is the universal C*-algebra generated by a copy of $A$ and an
isometry $S$ subject to the relations
  \izitem
  \zitem $S a = \a(a)S$,
  \zitem $S^*aS = \Tr(a)$,  and
  \zitem $1 = \sum_{i=1}^m u_iS S^* u_i^*$,
  \medskip\noindent for all $a$ in $A$.

  \proof Follows immediately from the Theorem above.
  \proofend

We will take the remainder of this section to develop a few technical
consequences of the fact that $E$ is of index-finite type, to be used
in later sections.  Some of these are interesting in their own.  The
first one is essentially the combined contents of Propositions
\lcite{8.2} and \lcite{8.3} of \cite{\Tower}:

  \state Proposition
  \label \AuxOrderNRedundancy
  If\/ $\{u_1,\ldots,u_m\}$ is a quasi-basis for $E$ and $n\in \N$ then
  \izitem
  \zitem $\sum_{i=1}^m\a^n(u_i) S^{n+1}S^{*n+1}\a^n(u_i^*) =
    S^n S^{*n}.$
  \zitem $\K_n\subseteq\K_{n+1}$.

  \proof
  We have
  $$
  \sum_{i=1}^m\a^n(u_i) S^{n+1}S^{*n+1}\a^n(u_i^*) =
  \sum_{i=1}^m S^n u_i SS^* u_i^* S^{*n} =
  S^n S^{*n},
  $$
  proving (i).  The second point then follows immediately from (i).
  \proofend

  It follows that the fixed point algebra for the gauge action, namely
$\U$, is the inductive limit of the $\K_n$.
  Recall from \scite{\Tower}{8.9} that there exists a conditional
expectation $G : \CP \to A$ such that
  $$
  G(aS^nS^{*m}b) =
  \delta_{nm}aI_n\inv b
  \for a,b\in A \for n,m\in\N,
  $$
  where $\delta$ is the Kronecker symbol and
  $$
  I_n =
  \ind(E) \a\big(\ind(E)\big) \ldots \a^{n-1}\big(\ind(E)\big).
  \eqno {(\seqnumbering)}
  \label \DefIn
  $$

  % To see that $I_n$ is indeed a positive element observe that,
  % because $\ind(E)$ is central, $\a^p\big(\ind(E)\big)$ and
  % $\a^q\big(\ind(E)\big)$ are commuting elements of $A$, for all
  % $p,q\in\N$.

  \state Proposition
  \label \MainOrderNRedundancy
  Let $\{u_1,\ldots,u_m\}$ be a quasi-basis for $E$ and put $Z =
\{1,\ldots,m\}$.  For each $n\in\N\setminus\{0\}$ and
  for each multi-index $i = (i_0,i_1,\ldots,i_{n-1}) \in Z^n$ let
  $$
  u_{(i)} = u_{i_0} \a(u_{i_1}) \a^2(u_{i_2}) \ldots
\a^{n-1}(u_{i_{n-1}}).
  $$
  Then
  \izitem
  \zitem $\ds\sum_{i\in Z^n} u_{(i)}S^n S^{*n} u_{(i)}^*=1.$
  \zitem If $A$ is commutative then $\ds\sum_{i\in Z^n}
u_{(i)}u_{(i)}^* =
  I_n$.

  \proof
  For $n=1$ the conclusion follows from \lcite{\MainRedundancy} and
the definition of $\ind(E)$.  So assume that $n>1$ and observe that
  $$
  \sum_{i\in Z^n} u_{(i)}S^n S^{*n} u_{(i)}^* =
  \sum_{i\in Z^{n-1}} \sum_{j=1}^m u_{(i)}\a^{n-1}(u_j)S^n
S^{*n}\a^{n-1}(u_j^*) u_{(i)}^* \={\AuxOrderNRedundancy.i}
  \sum_{i\in Z^{n-1}} u_{(i)}S^{n-1}S^{*n-1}u_{(i)}^* = 1,
  $$
  where the last step is by induction.  As for (ii) observe that
  $$
  \sum_{i\in Z^n} u_{(i)}u_{(i)}^* =
  \sum_{i\in Z^{n-1}} \sum_{j=1}^m
    u_{(i)}\a^{n-1}(u_j)\a^{n-1}(u_j^*) u_{(i)}^* \$=
  % \sum_{i\in Z^{n-1}}
  % u_{(i)} \a^{n-1}\(\sum_{j=1}^mu_ju_j^*\) u_{(i)}^* \$=
  \sum_{i\in Z^{n-1}}
    u_{(i)} \a^{n-1}\big(\ind(E)\big) u_{(i)}^* =
  \a^{n-1}\big(\ind(E)\big) \sum_{i\in Z^{n-1}}
    u_{(i)} u_{(i)}^*,
 $$
  and the conclusion again follows by induction.
  \proofend

  As a consequence we will see that the conditional expectation $G$
may be described by an algebraic expression on each $K_n$ when $A$ is
abelian.

  \state Corollary
  \label \ImplemantsGonKn
  Suppose that $E$ is of index-finite type and $A$ is a commutative
C*-algebra.  Then
  for each $n\in\N$ there exists a finite set
$\{v_1,\ldots,v_p\}\subseteq A$ such that
  \izitem
  \zitem $G(a) = \sum_{i=1}^p v_i a v_i^*$, for all $a\in\K_n$, and
  \zitem $\sum_{i=1}^p v_i v_i^* = 1$.

  \proof
  Let $\{u_{(i)}\}_{i\in Z^n}$ be as in \lcite{\MainOrderNRedundancy}
and
  for $i\in Z^n$ set $v_{(i)} = I_n^{-1/2} u_{(i)}$, where $I_n$ is
defined in \lcite{\DefIn}.  Then  for any
  $a,b\in A$ we have
  $$
  \sum_{i\in Z^n} v_{(i)} a S^n S^{*n}b v_{(i)}^* =
  a I_n^{-1/2}\(\sum_{i\in Z^n} u_{(i)} S^n S^{*n} u_{(i)}^*\)
I_n^{-1/2} b=
  a I_n\inv b = G(a S^n S^{*n}b).
  $$
  Finally observe that plugging $a=1$ in (i) gives (ii).
  \proofend

  Recall that $\E_n$ was defined to be the composition
  $\E_n=\a^n\circ\Tr^n$.  If we observe that $\Tr=\a\inv\circ E$ we may
write $\E_n$ as

  $$
  \E_n =
  \a^{n} \underbrace{(\a\inv E) \ldots (\a\inv E)}_{n\rm\;times} =
  \(\a^{n-1} E \a^{-(n-1)}\) \ldots \(\a^2 E \a^{-2}\) \(\a E
\a^{-1}\) E.
  $$
  Viewing each $\a^k E \a^{-k}$ as a conditional expectation from
$\a^k(A)$ to $\a^{k+1}(A)$ it is obvious that it is of index-finite
type.  By \scite{\Watatani}{1.7.1} $\E_n$ is of index-finite type as
well and hence by
  \scite{\Watatani}{2.1.5} there are constants $\lambda_n>0$ such that
  $\|\E_n(a^*a)\|^{1/2} \geq \lambda_n \|a\|$, for all $a$ in $A$.

  \state Lemma
  Suppose that $E$ is of index-finite type.  Then
  for each $n\in \N$ one has that $AS^n$ is closed in $\CP$ so that
$\M_n = AS^n$ (without closure).

  \proof For each $a\in A$ we have
  $$
  \|a S^n\| =
  \|S^{*n}a^*a S^n\|^{1/2} =
  \|\Tr^n(a^*a)\|^{1/2} =
  \|\a^n(\Tr^n(a^*a))\|^{1/2} =
  \|\E_n(a^*a)\|^{1/2} \geq
  \lambda_n \|a\|.
  $$
  Therefore the map $a\in A \mapsto aS^n\in \M_n$ is a Banach space
isomorphism onto its range which in turn is a complete normed space,
hence closed.
  \proofend

  \state Proposition
  \label \PropertiesOfFinite
  Suppose that $E$ is of index-finite type and
  let $n\in\N$.  Then any element $k\in \K_n$ may be written as a
finite sum of the form $\sum_{i=1}^m a_iS^nS^{*n}b_i^*$, where $m\in
\N$ and $a_i,b_i\in A$.
  If $k\geq0$ we may take $b_i=a_i$.

  \proof
  Let $\{u_{(i)}\}_{i\in Z^n}$ be as in \lcite{\MainOrderNRedundancy}.
Then for each $i\in Z^n$ we have by \lcite{\MisModuleOverK.i} that $k
u_{(i)} S^n\in \M_n = AS^n$.  So there is $a_i\in A$ such that
  $k u_{(i)} S^n = a_i S^n$ and by \lcite{\MainOrderNRedundancy.i}
  $$
  k =
  k \(\sum_{i\in Z^n} u_{(i)}S^nS^{*n}u_{(i)}^*\) =
  \sum_{i\in Z^n} a_iS^nS^{*n}u_{(i)}^*.
  $$
  If $k \geq 0$ write $k = l^*l$, with
$l=\sum_{i=1}^mx_iS^nS^{*n}y_i$, and $x_i,y_i\in A$.
  Then
  $$
  k =
  \sum_{i,j=1}^m x_iS^nS^{*n}y_iy_j^*S^nS^{*n}x_j^* =
  \sum_{i,j=1}^m x_i\E_n(y_iy_j^*)S^nS^{*n}x_j^*.
  $$
  Since $\big(\E_n(y_iy_j^*)\big)_{i,j}$ is a positive $m\times m$
matrix over $\a^n(A)$ by
  \scite{\Takesaki}{IV.3.4}
  there exists another such matrix, say $\big(c_{i,j}\big)_{i,j}$ such
that
  $\E_n(y_iy_j^*) = \sum_{k=1}^m c_{ik}c_{jk}^*$
  for all $i$ and $j$.  Therefore
  $$
  k =
  \sum_{i,j,k=1}^m x_i c_{ik}c_{jk}^* S^nS^{*n}x_j^* =
  \sum_{i,j,k=1}^m x_i c_{ik}S^nS^{*n}c_{jk}^* x_j^* \$=
  \sum_{k=1}^m \(\sum_{i=1}^m x_i c_{ik}\) S^nS^{*n} \(\sum_{j=1}^m
c_{jk}^* x_j^*\) =
  \sum_{k=1}^m a_k S^nS^{*n} a_k^*,
  $$
  where $a_k = \sum_{i=1}^m x_i c_{ik}$.
  \proofend

  \state Proposition
  \label \GfaithfulOnKn
  If $E$ is of index-finite type the restriction of the conditional
expectation $G$ above to each $\K_n$ is faithful.

  \proof
  Let $k\geq0$ in $\K_n$ be such that $G(k)=0$.  Write
  $k = \sum_{i=1}^m a_iS^nS^{*n}a_i^*$ as in
\lcite{\PropertiesOfFinite}.  Then
  $$
  0 =
  G\(\sum_{i=1}^m a_iS^nS^{*n}a_i^*\) =
  \sum_{i=1}^m a_i I_n a_i^*,
  $$
  which implies that $a_i I_n a_i^*=0$, for all $i$, and hence also
that $a_i I_n^{1/2}=0$.  Since $I_n$ is invertible we have that $a_i=0$ and
so $k=0$.
  \proofend

  We will now see that $G$ is also faithful on $\CP$ in the following
weak sense.

  \state Proposition
  \label \GfaithfulOnU
  Suppose that $E$ is of index-finite type and
  let $a\geq 0$ in $\CP$ be such that $G(xax^*)=0$ for all $x\in \U$.
Then $a=0$.

  \proof
  Using the conditional expectation $F$ of \lcite{\CondExpF} we have
  $$
  0 = G(xax^*) = G(F(xax^*)) = G(xF(a)x^*),
  $$
  for all $x\in \U$.
  Observe that this allows us to suppose that $a\in \U$.  In fact,
should the result be proved in this special case, we would have that
$F(a)=0$ in which case $a=0$ since $F$ is faithful.  We therefore
suppose that $a\in \U$.

  By the polarization identity we have that $G(xay^*)=0$ for all
$x,y\in \U$.  Supposing by way of contradiction that $a\neq0$ observe
that the closed two sided ideal $J$ of $\U$ given by
  $$
  J = \{ b \in \U : G(xby^*)=0,\ \forall x,y\in \U\}.
  $$
  is nontrivial.
  Since $\U = \overline{\bigcup_{n\in \N} \K_n}$ it follows from
\scite{\alnr}{1.3} (see also \scite{\Bratteli}{3.1}) that there exists
a nonzero element $b \in J\cap \K_n$ for some $n$.  But since
$G(b^*b)=0$ we have by \lcite{\GfaithfulOnKn} that $b=0$, which is a
contradiction.
  % Note to myself: I would like to prove that $G$ is faithful on the
  % inductive limit \break but I can't!
  \proofend

  \section {Covering maps}
  \label \CommutativeSection
  From this section on we will let $X$ be a compact topological space
and $T:X\to X$ be a continuous surjective map.  We will then consider
the C*-algebra $A=C(X)$ equipped with the endomorphism $\a$ given by
  $$
  \a: f\in C(X) \mapsto f \circ T \in C(X).
  \eqno {(\seqnumbering)}
  \label \DefEndoForCovering
  $$

Another important hypothesis that we will assume  is that
$T$ is a covering map.  As we will see this will allow us to
define a somewhat canonical transfer operator.

  \state Proposition
  If $T:X\to X$ is a covering map then for each $f\in A$ the function
  $\L(f)$ given by
  $$
  \L(f)\calcat x = \sum_{t\in T\inv(x)} f(t)
  \for x\in X,
  $$
  is continuous on $X$.

  \proof
  Given $x_0\in X$ let $W$ be an open neighborhood of $x_0$ such that
$T\inv(W)$ is the disjoint union of open sets $U_1,\ldots,U_n$ such
that $T$ is a homeomorphism from each $U_i$ onto $W$.  For each
$i=1,\ldots,n$ let $\phi_i: W \to U_i$ be the inverse of the
restriction of $T$ to $U_i$.  Then for all $x\in W$ one has that
$T\inv(x) = \{\phi_1(x),\ldots,\phi_n(x)\}$ and hence
  $$
  \L(f)\calcat x = \sum_{i=1}^n f(\phi_i(x))
  $$
  so that $\L(f)$ is seen to be continuous at $x_0$.
  \proofend

  Setting
  $$
  \Tr(f) = \L(1)\inv \L(f)
  \and
  E = \a \circ \Tr
  \eqno {(\seqnumbering)}
  \label \DefTransForCovering
  $$
  it is easy to see that $\Tr$ is a transfer operator for $\a$,
  that $E$ is a conditional expectation from $A$ onto the range of
$\a$, and that $\Tr=\a\inv\circ E$.  Adopting the notation $\I =
\a(\L(1))$ it is easy to see that for all $x\in X$
  $$
  \I(x) = \# \big\{t\in X : T(t)=T(x)\big\},
  \eqno {(\seqnumbering)}
  \label \ThisIsIndex
  $$
  and
  $$
  E(f)\calcat x =
  {1 \over \I(x)}
  \sum_{\buildrel {\scriptstyle t\in X} \over {T(t)=T(x)}} f(t).
  \eqno {(\seqnumbering)}
  \label \HereIsMyExp
  $$

  The next result is taken from \scite{\Tower}{Section 11} and is
reproduced here for the convenience of the reader:

  \state Proposition
  \label \CoveringIsFinite
  Let $\{V_i\}_{i=1}^m$ be a finite open covering of $X$ such that
the restriction of\/ $T$ to each  $V_i$ is one-to-one
and let $\{v_i\}_{i=1}^m$ be
a partition of unit subordinate to that covering.  Setting $u_i=(\I
v_i)^{1/2}$ one has that $\{u_i\}_{i=1}^m$ is a quasi basis for $E$
and hence $E$ is of index-finite type.  Moreover $\ind(E) = \I$.

  \proof
  Observe that for all $f\in A$ and $x\in X$ one has that
  $$
  \sum_{i=1}^m u_i E(u_if)\calcat x =
  \sum_{i=1}^m u_i(x) {1 \over \I(x)}
    \sum_{\buildrel {\scriptstyle t\in X} \over {T(t)=T(x)}}
u_i(t)f(t) \$=
  \sum_{i=1}^m u_i(x) {1 \over \I(x)} u_i(x)f(x)=
  \sum_{i=1}^m v_i(x) f(x) = f(x).
  $$
  Therefore $\{u_1,\ldots,u_m\}$ is a quasi-basis for $E$ and
  $$
  \ind(E) =
  \sum_{i=1}^m u_i^2 =
  \sum_{i=1}^m \I v_i = \I.
  \proofend
  $$

In the final result of this section we denote by $\supp(f)$ the support
of the function $f$, namely the closure of the set $\{x : f(x) \neq
0\}$.

  \state Lemma
  \label \SupportLemma
  Let $\Tr$ be any transfer operator for $\a$ (such as the one defined
in \lcite{\DefTransForCovering} but not necessarily).  Then
for $f\in C(X)$ one has
  \izitem
  \zitem $\supp\big(\a(f)\big) \subseteq T\inv\big(\supp(f)\big)$,
and.
  \zitem $\supp\big(\Tr(f)\big) \subseteq T\big(\supp(f)\big)$.

  \proof
  The first point is trivial. As for (ii)
  let $x$ be such that $\Tr(f)\calcat x\neq 0$.  We claim that $x\in
T(\supp(f))$.  Otherwise let $g$ be a continuous function on $X$ which
vanishes on $T(\supp(f))$ and such that $g(x)=1$.  Then $f\a(g)=0$
because for all $y\in X$
  $$
  f(y)\, \a(g)(y) =
  f(y)\, g(T(y)),
  $$
  and if $f(y)\neq 0$ then $y\in \supp(f)$ in which case $g(T(y))=0$.
We then have
  $$
  \Tr(f)\calcat x =
  \Tr(f)\calcat x\, g(x) =
  \Tr(f)g\calcat x =
  \Tr\big(f\a(g)\big)\calcat x =0,
  $$
  contradicting our assumption.  This proves our claim that $x\in
T(\supp(f))$.
  Since $T(\supp(f))$ is closed the conclusion follows.
  \proofend

  \section{The algebra $\CA(X,T)$}
  As before let $X$ be a compact topological space and let $T:X\to X$
be a covering map.  We will suppose, as in \cite{\Deaconu}, that each
point of $X$ has exactly $p$ inverse images under $T$.

Motivated by earlier work of Renault on the Cuntz groupoid
\cite{\Renault} and on Arzumanian and Vershik's paper
\cite{\ArzVershikTwo}, Deaconu \cite{\Deaconu} has considered the set
  $$
  \G = \{(x,n,y)\in X\times\Z \times X: \exists\, k,l \in \N,\ n=l-k,\
T^k(x) = T^l(y)\}
  $$
  with the following groupoid structure: the product of $(x,n,y)$ and
$(w,m,z)$ is defined iff $y=w$ in which case it is given by
$(x,n+m,z)$.  The inverse of $(x,n,y)$ is set to be $(y,-n,x)$.
Deaconu has proven that $\G$ admits the structure of an $r$-discrete
locally compact groupoid with Haar system given by counting measures
(see \scite{\Deaconu}{Theorem 1}).  We will denote the corresponding
groupoid C*-algebra \cite{\Renault} by $\CA(X,T)$.

We will view $C(X)$ as a subalgebra of $\CA(X,T)$ by identifying a
function $f\in C(X)$ with the function (also denoted) $f$ on $\G$
given by
  \def\bool#1{\big[#1\big]\ }
  $$
  f(x,n,y) = \bool{n=0 \ \wedge\ x = y}
  f(x),
  $$
  where the brackets correspond to the obvious boolean valued
function.  Let $v$ be element of $\CA(X,T)$ described shortly before
the statement of Theorem 2 in \cite{\Deaconu} by
  $$
  v(x,n,y) = \bool{n=-1 \ \wedge \ y = T(x)} / \sqrt p
  $$
  where $p$ is the index of the covering.

  \state Theorem 
  \label \CoveringIsomCrossProd
  Let $X$ be a compact topological space and $T:X\to X$ be a covering map.
  Let $\a$ be the endomorphism of\/ $C(X)$ given by
\lcite{\DefEndoForCovering} and let $\Tr$ be the transfer operator
defined in \lcite{\DefTransForCovering}.  Then there exists a
*-isomorphism $\phi$ from $C(X)\crpr$ onto $\CA(X,T)$ which is the
identity on $C(X)$ and such that $\phi(S)=v$.

  \proof
  By \cite{\Deaconu} we have that
  $$
  (vv^*)(x,n,y) = \bool{n=0 \ \wedge\ T(x) = T(y)}/p.
  $$
  Let $v_i$, $u_i$, and $\I$ be as in \lcite{\CoveringIsFinite}
observing that in the present case $\I$ is the constant function equal
to $p$.  Since $T$ is injective when restricted to the set where each
$u_i$ does not vanish we have
  $$
  (u_ivv^* u_i)(x,n,y) =
  u_i(x)\ (vv^*)(x,n,y)\ u_i(y) =
  \bool{n=0 \ \wedge\ T(x) = T(y)}
  u_i(x)
  u_i(y)/p \$=
  \bool{n=0 \ \wedge\ x = y} u_i(x)^2 /p =
  \bool{n=0 \ \wedge\ x = y} v_i(x).
  $$
  Since the $v_i$ form a partition of unit we have that
  $$
  \sum_{i=1}^m u_ivv^* u_i = 1.
  $$
  By direct computation one can prove that $vf = \a(f)v$ and $v^*f v =
\Tr(f)$ for each $f\in C(X)$ so by \lcite{\ShorterRelations} there
exists a *-homomorphism
  $$
  \phi : C(X)\crpr \to \CA(X,T)
  $$
  such that $\phi(f)=f$, for all $f\in C(X)$ and $\phi(S)=v$.  It is
clear that $\phi$ is injective on $C(X)$ and covariant for the gauge
action on $C(X)\crpr$ and the circle action on $\CA(X,T)$ described in
\cite{\Deaconu} in terms of the cocycle
  $$
  c(x,n,y)= -n.
  $$
  That $\phi$ is injective then follows from
\lcite{\GaugeInvariantUniqueness}.  On the other hand we have that
$C(X)\cup \{v\}$ generates $\CA(X,T)$ as a C*-algebra and hence $\phi$
is surjective as well.
  \proofend

  The description of the $\CA(X,T)$ in terms of generators and
relations therefore becomes:

  \state Theorem
  \label \ShorterRelationsForCovering
  Let $X$ be a compact topological space and let $T:X\to X$ be a
covering map.  Then $\CA(X,T)$ is the universal C*-algebra generated
by a copy of\/ $C(X)$ and an isometry $S$ subject to the relations
  \izitem
  \zitem $S f = \a(f)S$,
  \zitem $S^*fS = \Tr(f)$,  and
  \zitem $1 = \sum_{i=1}^m u_iS S^* u_i^*$,
  \medskip\noindent for all $f$ in $C(X)$,
  where $\Tr$ is given by \lcite{\DefTransForCovering} and 
the $u_i$ are as in \lcite{\CoveringIsFinite}.

  \proof Follows immediately from the result above and
\lcite{\ShorterRelations}.
  \proofend

  Before closing this section we should remark that an example of the
situation treated above is given by the Markov one-sided subshift
associated to a zero-one matrix $A$ with no zero rows, in which case
$\CCP$ is isomorphic to the Cuntz-Krieger algebra ${\cal O}_A$ by
\scite{\Endo}{Theorem 6.2}.

  \section{Topological Freeness}
  In this section we will obtain another result about faithfulness of
representations.  This time, instead of requiring gauge-covariance as
in \lcite{\GaugeInvariantUniqueness}, we will make assumptions on the
dynamical properties of the transformation $T$.

  We keep the standing assumptions of section
\lcite{\CommutativeSection}, namely that $T$ is a covering map of the
compact space $X$ and the endomorphism $\a$ of $C(X)$ is given by
\lcite{\DefEndoForCovering}.  The transfer operator $\Tr$ will be
given, as before, by \lcite{\DefTransForCovering}.

We begin by describing the main hypothesis to be used below (see also
  \scite{\Tom}{2.1},
  \cite{\AS},
  \cite{\Deaconu}, and
  \scite{\ELQ}{2.1}).

  \definition We will say that the dynamical system $(X,T)$ is
\stress{topologically free} if for every pair of nonnegative integers
$(n,m)$ with $n\neq m$ one has that the set
  $$
  \{ x \in X : T^n(x) = T^m(x) \}
  $$
  has empty interior.

  The next result is motivated by Lemma \lcite{2.3} in \cite{\ELQ}.

  \state Lemma
  \label \GrandeH
  Let $x_0\in X$ and $n,m\in\N$ be such that $T^n(x_0) \neq T^m(x_0)$.
Then there exists $h\in C(X)$ with
  \izitem
  \zitem $0 \leq h \leq 1$,
  \zitem $h(x_0)=1$,
  \zitem $h S^n S^{*m} h = 0$.

  \proof
  We begin by claiming that there exists an open set $U\subseteq X$
with $x_0\in U$ such that $U \cap T^{-n}(T^m(U)) = \emptyset$.  In
fact let $A,B\subseteq X$ be disjoint open sets such that $T^n(x_0)
\in A$ and $T^m(x_0) \in B$ and set $U = T^{-n}(A) \cap T^{-m}(B)$.
Then obviously $x_0\in U$ and $T^n(U) \cap T^m(U) =\emptyset$.  In
addition
  $$
  \emptyset =
  T^{-n}\big(T^n(U) \cap T^m(U)\big) =
  T^{-n}\big(T^n(U)\big) \cap T^{-n}\big(T^m(U)\big) \supseteq
  U \cap T^{-n}\big(T^m(U)\big),
  $$
  thus proving our claim.
  Pick $h\in C(X)$ satisfying (i) and (ii) above and such that
$\supp(h)\subseteq U$.  Then
  $$
  \|h S^n S^{*m} h\|^2 =
  \|h S^n S^{*m} h^2 S^m S^{*n} h\| =
  \|h \a^n(\Tr^m(h^2)) S^n S^{*n} h\|.
  $$
  So it suffices to prove that $h \a^n(\Tr^m(h^2))=0$.  For this
purpose observe that by \scite{\Endo}{3.2} we have that $\Tr^m$ is a
transfer operator for $\a^m$ and hence by \lcite{\SupportLemma} we
have
  $$
  \supp\Big(\a^n\big(\Tr^m(h^2)\big)\Big) \subseteq
  T^{-n}\Big(\supp\big(\Tr^m(h^2)\big)\Big) \subseteq
  T^{-n}\Big(T^m\big(\supp(h^2)\big)\Big) \subseteq
  T^{-n}(T^m(U)).
  $$
  Since the latter set is disjoint from $U$, and hence also from
$\supp(h)$, we have that $h \a^n(\Tr^m(h^2))$ is indeed zero.
  \proofend

  We now come to the main result of this section (for similar results
see
  \scite{\Stacey}{2.1},
  \scite{\BKR}{2.1},
  \scite{\alnr}{1.2},
  \scite{\quasilat}{3.7},
  \scite{\ELQ}{2.6}).

  \state Theorem
  \label \MainResultTopFree
  Let $T$ be a topologically free covering map of the compact space
$X$, let $\a$ be the endomorphism of $C(X)$ given by
\lcite{\DefEndoForCovering}, and $\Tr$ be the transfer operator given
by \lcite{\DefTransForCovering}.
  Then any nontrivial ideal of $\CCP$ must have a nontrivial
intersection with $C(X)$.
  Given a C*-algebra $B$ and a *-homomorphism $\pi: \CCP \to B$ which
is faithful on $C(X)$ one has that $\pi$ is faithful on all of $\CCP$.

  \proof
  Let $\pi: \CCP \to B$ be faithful on $C(X)$.
  We begin by claiming that
  $$
  \|G(a)\| \leq \|\pi(a)\|
  \eqno {(\dagger)}
  $$
  for all $a\geq0$ in $\CCP$.  In order to prove this observe that the
set of nonnegative elements of the form
  $$
  a = \sum_{i=1}^t a_i S^{n_i}S^{*m_i} b_i,
  \eqno {(\ddagger)}
  $$
  where $n_i,m_i\in\N$ and $a_i,b_i\in C(X)$, is dense in the positive
cone of $\CCP$.  In fact, given any $a\geq0$ in $\CCP$ one may write
$a = b^*b$ and then approximate $b$ by elements of the above form by
\scite{\Tower}{2.3}.  Using \scite{\Tower}{2.2} we have that $b^*b$
again has the above form and clearly $b^*b$ can be made to be as close
as necessary to $a$.  In order to prove $(\dagger)$ we may therefore
assume that $a$ has the form of $(\ddagger)$.
  Let us also suppose that the sum in $(\ddagger)$ is arranged in such
a way that $n_i=m_i$ for all $i=1,\ldots,s$, while $n_i\neq m_i$ for
$i=s+1,\ldots,t$.

For each $i=s+1,\ldots,t$ we have by hypothesis that $\{x \in X :
T^{n_i}(x) \neq T^{m_i}(x)\}$ is an open dense set and hence so is
their intersection.  Consequently, fixing a real number $\rho$ with
$0<\rho<1$, there exists $x_0\in X$ such that $G(a)\calcat{x_0} >
\rho\|G(a)\|$ and $T^{n_i}(x_0) \neq T^{m_i}(x_0)$ for all
$i=s+1,\ldots,t$.

Using \lcite{\GrandeH} choose for each $i=s+1,\ldots,t$ an $h_i\in
C(X)$ with $0 \leq h_i \leq 1$, $h_i(x_0)=1$, and $h_i S^{n_i}
S^{*m_i} h_i = 0$.  Setting $h = h_{s+1}\cdots h_t$ we therefore have
that
  $$
  h a h =
  \sum_{i=1}^t a_i h S^{n_i}S^{*m_i}h b_i =
  \sum_{i=1}^s a_i h S^{n_i}S^{*n_i}h b_i =
  h F(a) h,
  $$
  where $F$ is given by \lcite{\CondExpF}.
  Since $F(a)\in \K_n$, where $n = \max\{n_1,\ldots,n_s\}$, by
\lcite{\AuxOrderNRedundancy.ii} we may choose by
\lcite{\ImplemantsGonKn} a finite set $\{v_1,\ldots,v_p\}\subseteq
C(X)$ such that $\sum_{i=1}^p v_i v_i^* = 1$ and
  $$
  G(a) = G(F(a)) = \sum_{i=1}^p v_i F(a) v_i^*.
  $$
  It follows that
  $$
  \| hG(a)h \| =
  \| \pi\big(hG(a)h\big) \| =
  \left\| \sum_{i=1}^p \pi\big(h v_i F(a)v_i^*h\big)\right\| =
  \left\| \sum_{i=1}^p \pi(v_i hah v_i^*)\right\| \leq
  \|\pi(a)\|.
  $$
  Since
  $$
  \| hG(a)h \| \geq
  h(x_0)G(a)\calcat{x_0} h(x_0) =
  G(a)\calcat{x_0} >
  \rho\|G(a)\|
  $$
  and $\rho$ is arbitrary we conclude that $\|G(a)\| \leq \|\pi(a)\|$
as claimed.

  Let $a\in\CCP$ be such that $\pi(a)=0$.  Then for any $b\in\CCP$ we
have
  $$
  \|G(b^*a^*ab)\| \leq \|\pi(b^*a^*ab)\| =0,
  $$
  and hence $G(b^*a^*ab) =0$. Applying
  \lcite{\GfaithfulOnU} we conclude that $a^*a=0$ and hence also that
$a=0$.

  The assertion about ideals in the statement is proved as follows.
Given an ideal $J$ of $\CCP$ which is not the trivial ideal let $\pi$
be the quotient map.  Then $\pi$ cannot be faithful on $C(X)$ or else
it would be faithful on all of $\CCP$ by what we have already proved.
It follows that $J\cap C(X) \neq \{0\}$.
  \proofend

  \section{Simplicity}
  As before we fix a compact space $X$, a covering map $T$ of $X$, and
we let $\a$ and $\Tr$ be given respectively by
\lcite{\DefEndoForCovering} and \lcite{\DefTransForCovering}.
  Our main goal here will be to find a necessary and sufficient
condition for $\CCP$ to be simple.

  Recall that two points $x,y\in X$ are said to be
\stress{trajectory-equivalent} (see e.g.~\cite{\ArzVershik}) when
there are $n,m\in\N$ such that $T^n(x) = T^m(y)$.  We will denote this
by $x\sim y$.
  A subset $Y\subseteq X$ is said to be invariant if $x\sim y \in Y$
implies that $x\in Y$.  It is easy to see that $Y$ is invariant if and
only if\/ $T\inv(Y)=Y$.
  We will say that $T$ is \stress{irreducible} when there is
no closed (equivalently open) invariant set other than $\emptyset$ and
$X$.  Notice that irreducibility is weaker than the
condition of minimality defined in \cite{\Deaconu}.

  \state Proposition
  \label \IrreducibleImplyTopFree
  If\/ $T$ is irreducible then either
  \izitem
  \zitem $T$ is topologically free, or
  \zitem $X$ is finite and $T$ is a cyclic permutation of $X$.

  \proof
  Suppose that $T$ is irreducible and not topologically free.  Then there is a
nonempty open set $U\subseteq X$ and $n,m\in\N$ with $n\neq m$ such
that $T^n(u)=T^m(u)$ for all $u$ in $U$.
  Clearly $\bigcup_{k,l\in\N} T^{-k}(T^l(U))$ is an open invariant set
so, being nonempty, it must coincide with $X$.  In this case
$\big\{T^{-k}(T^l(U))\big\}_{k,l\in\N}$ is an open cover of $X$ and
hence admits a finite subcover, say
  $\big\{T^{-k_i}(T^{l_i}(U))\big\}_{i=1,\ldots,p}$.
  Given any $x\in X$ choose $i$ such that $x\in T^{-k_i}(T^{l_i}(U))$.
Therefore there exists $u\in U$ such that $T^{k_i}(x) = T^{l_i}(u)$.
It follows that
  $$
  T^{n+k_i}(x) =
  T^{n}(T^{l_i}(u)) =
  T^{l_i}(T^{n}(u)) =
  T^{l_i}(T^{m}(u)) =
  T^{m}(T^{l_i}(u)) =
  T^{m+k_i}(x).
  $$
  Setting $k = \max\{k_1,\ldots,k_p\}$ we therefore see that
  $$
  T^{n+k}(x) = T^{m+k}(x)
  \for x\in X.
  $$
  Assuming without loss of generality that $n>m$, and setting $r=m+k$
and $s = n-m$, this translates into $T^{s+r}(x) = T^{r}(x)$.  Since
$T^r$ is a surjective map we conclude that $T^s=1$.  The reader may
now easily prove that we must be in the situation described by (ii).
  \proofend

  The next result improves on Deaconu's characterization of simplicity
\cite{\Deaconu}.  It is also a generalization of the corresponding
result for crossed-products by automorphisms (see \cite{\EffrosHahn},
\cite{\Power}, \cite{\Tom}).

  \state Theorem 
  \label \SimplicityResult
  Let $T$ be a covering map of an infinite compact space $X$ to
itself, let $\a$ be the endomorphism of $C(X)$ given by
\lcite{\DefEndoForCovering}, and let $\Tr$ be the transfer operator
given by \lcite{\DefTransForCovering}.  Then $\CCP$ is simple if and
only if $T$ is irreducible.

  \proof
  Supposing that $T$ is irreducible
  let $J$ be an ideal in $\CCP$.  Obviously $J\cap C(X)$ is an ideal
in $C(X)$ and hence coincides with $C_0(U)$ for some open set
$U\subseteq X$.  We claim that $U$ is invariant.  For this we need to
prove that $x\sim y \in U$ implies that $x\in U$.

  Let $f$ be a nonnegative function in $C_0(U)$ such that $f(y)=1$.
Given that $x\sim y$, let $n,m\in\N$ be such that $T^n(x) = T^m(y).$
  Taking $\{u_{(i)}\}_{i\in Z^n}$ as in \lcite{\MainOrderNRedundancy}
we therefore have that
  $$
  J \ni
  \sum_{i\in Z^n} u_{(i)} S^nS^{*m}fS^mS^{*n} u_{(i)}^* =
  \sum_{i\in Z^n} u_{(i)} \a^n\big(\Tr^m(f)\big)S^nS^{*n} u_{(i)}^*
\$=
  \a^n\big(\Tr^m(f)\big) \sum_{i\in Z^n} u_{(i)} S^nS^{*n} u_{(i)}^* =
  \a^n\big(\Tr^m(f)\big),
  $$
  by \lcite{\MainOrderNRedundancy.i}.
  It follows that $\a^n\big(\Tr^m(f)\big)\in C_0(U)$.  On the other
hand observe that
  $$
  \a^n\big(\Tr^m(f)\big)\calcat x =
  \Tr^m(f)\calcat {T^n(x)} =
  \Tr^m(f)\calcat {T^m(y)} =
  \a^m\big(\Tr^m(f)\big)\calcat y =
  \E_m(f)\calcat y \geq
  c f(y) > 0,
  $$
  where the constant $c$ appearing above is obtained from
\scite{\Watatani}{2.1.5}.  One must then have that $x\in U$, which
concludes the proof of the claim that $U$ is invariant.
  Since $T$ is irreducible we have that $U$ is either empty or equal
to $X$.  In case $U=\emptyset$ we have that $J\cap C(X)=\{0\}$, but
since $T$ is topologically free by \lcite{\IrreducibleImplyTopFree} we
have that $J=\{0\}$ by \lcite{\MainResultTopFree}.  On the other hand
if $U=X$ we have that $1\in C_0(U) \subseteq J$ in which case
$J=\CCP$.

   Let us now suppose that $T$ is not irreducible, that is, that there
exists a nontrivial open invariant subset $U\subseteq X$.  Let $I =
C_0(U)$.  It is then clear that $\a(I)\subseteq I$ and
$\Tr(I)\subseteq I$.  The conclusion then follows from the following
result which actually holds under more general hypothesis.
  \proofend

  \state Proposition
  Let $\a$ be an injective *-endomorphism of a unital C*-algebra $A$
with $\a(1)=1$ and let $\Tr$ be a transfer operator for $\a$ such that
$\Tr(1)=1$.  If there exists a closed two-sided ideal $I$ of $A$ such
that $\a(I)\subseteq I$ and $\Tr(I)\subseteq I$ then there exists a
closed two-sided ideal $J$ of $\CP$ such that $J\cap A = I$.

  \proof
  Consider the operators $\tilde\a$ and $\tilde\Tr$ on $A/I$ obtained
by passing $\a$ and $\Tr$ to the quotient, respectively.  It is then
obvious that $\tilde\a$ is an endomorphism of $A/I$ and $\tilde\Tr$
is a transfer operator for $\tilde\a$.  Denote by $\tilde S$ the
standard isometry of $(A/I) \crossproduct_{\tilde\a,\tilde\Tr} {\bf
N}$ and view the quotient map of $A$ modulo $I$ as a *-homomorphism
  $$
  q : A \to A/I \subseteq
  (A/I) \crossproduct_{\tilde\a,\tilde\Tr} {\bf N}.
  $$
  It is elementary to check that
  $$
  \tilde S q(a) = q(\a(a)) \tilde S
  \and
  \tilde S^* q(a) \tilde S = q(\Tr(a)),
  $$
  for all $a\in A$,
  which implies that there exists a *-homomorphism
  $$
  \phi: \TCP \to (A/I) \crossproduct_{\tilde\a,\tilde\Tr} {\bf N}
  $$
  such that $\phi (a) = q(a)$, for all $a\in A$, and $\phi(\tS) =
\tilde S$.  If $(a,k)\in A \times \overline{A \tS\tS^*A}$ is a
redundancy it is easy to see that $\phi(a) = \phi(k)$ and hence that
$\phi$ drops to the quotient providing a *-homomorphism
  $$
  \psi: \CP \to (A/I) \crossproduct_{\tilde\a,\tilde\Tr} {\bf N}
  $$
  which coincides with $\phi$, and hence also with $q$, on $A$.  The
kernel of $\psi$ is then the desired ideal.
  \proofend

  \endgroup


\references

\bibitem{\alnr}
  {S. Adji, M. Laca, M. Nilsen, and I. Raeburn}
  {Crossed products by semigroups of endomorphisms and the Toeplitz
algebras of ordered groups}
  {{\it Proc. Amer. Math. Soc.}, {\bf 122} (1994), 1133--1141}

\bibitem{\AS}
  {R. J. Archbold and J. S. Spielberg}
  {Topologically free actions and ideals in discrete C*-dynamical
systems}
  {{\it Proc. Edinb. Math. Soc., II. Ser. 37} (1994), 119--124}

\bibitem{\ArzVershik}
  {V. A. Arzumanian and A. M. Vershik}
  {Factor representations of the crossed product of a commutative
C*-algebra and a semigroup of its endomorphisms}
  {{\it Dokl. Akad. Nauk. SSSR}, {\bf 238} (1978), 513--517.
  Translated in {\it Soviet Math.~Dokl.}, {\bf 19} (1978), No.~1}

\bibitem{\ArzVershikTwo}
  {V. A. Arzumanian and A. M. Vershik}
  {Star algebras associated with endomorphisms}
  {Operator Algebras and Group Representations, Proc. Int. Conf.,
Neptun/Rom. 1980, Vol. I, Monogr. Stud. Math. 17, 17-27 (1984)}


\bibitem{\Bourbaki}
  {N. Bourbaki}
  {\'El\'ements de math\'ematique. XXV. Premi\`ere partie. Livre VI:
Int\'egration. Chapitre 6: Int\'egration vectorielle}
  {Hermann, 1959}


\bibitem{\BKR}
  {S. Boyd, N. Keswani, and I. Raeburn}
  {Faithful representations of crossed products by endomorphisms}
  {{\it Proc. Am. Math. Soc.}, {\bf 118} (1993), 427--436}

\bibitem{\Bratteli}
  {O. Bratteli}
  {Inductive limits of finite dimensional $C^*$-algebras}
  {{\it Trans. Amer. Math. Soc.}, {\bf 171} (1972), 195--234}

\bibitem{\Deaconu}
  {V. Deaconu}
  {Groupoids associated with endomorphisms}
  {{\it Trans. Amer. Math. Soc.}, {\bf 347} (1995), 1779--1786}

\bibitem{\EffrosHahn}
  {E. G. Effros and F. Hahn}
  {Locally compact transformation groups and C*-algebras}
  {{\it Mem. Am. Math. Soc.}, {\bf 75} (1967), 92 p}

\bibitem{\newpim}
  {R. Exel}
  {Circle Actions on C*-Algebras, Partial Automorphisms and a
Generalized Pimsner--Voiculescu Exact Sequence}
  {{\it J. Funct. Analysis}, {\bf 122} (1994), 361--401,
[arXiv:funct-an/9211001]}

\bibitem{\Endo}
  {R. Exel}
  {A New Look at The Crossed-Product of a C*-algebra by an
Endomorphism}
  {preprint, Universidade Federal de Santa Catarina, 2000,
[arXiv:math.OA/0012084]}

\bibitem{\Tower}
  {R. Exel}
  {Crossed-Products by Finite Index Endomorphisms and KMS states}
  {preprint, Universidade Federal de Santa Catarina, 2001,
[arXiv:math.OA/0105195]}

\bibitem{\RPF}
  {R. Exel}
  {KMS states for generalized gauge actions on Cuntz-Krieger algebras
(An application of the Ruelle-Perron-Frobenius Theorem)}
  {preprint, Universidade Federal de Santa Catarina, 2001,
[arXiv:math.OA/0110183]}

\bibitem{\ELQ}
  {R. Exel, M. Laca, and J. Quigg}
  {Partial Dynamical Systems and C*-Algebras generated by Partial
Isometries}
  {{\it J. Operator Theory}, to appear, [arXiv:funct-an/9712007]}

\bibitem{\anHuef}
  {A. an Huef and I. Raeburn}
  {The ideal structure of Cuntz-Krieger algebras}
  {{\it Ergodic Theory Dyn. Syst.}, {\bf 17} (1997), 611--624}

\bibitem{\quasilat}
  {M. Laca and I. Raeburn}
  {Semigroup crossed products and Toeplitz algebras of nonabelian
groups}
  {{\it J. Funct. Analysis}, {\bf 139} (1996), 415--440}

  % \bibitem{\Lance}
  % {E. C. Lance}
  % {Hilbert C*-modules. A toolkit for operator algebraists}
  % {London Mathematical Society Lecture Note Series vol.~210,
  % Cambridge Univ. Press, 1995}

\bibitem{\Murphy}
  {G. J. Murphy}
  {Crossed products of C*-algebras by endomorphisms}
  {\sl Integral Equations Oper. Theory \bf 24 \rm (1996), 298--319}

\bibitem{\VN}
  {F. J. Murray and J. von Neumann}
  {On rings of operators}
  {{\it Ann. of Math.}, {\bf 37} (1936), 116--229}

\bibitem{\Ped}
  {G. K. Pedersen}
  {$C^*$-Algebras and their automorphism groups}
  {Acad. Press, 1979}

\bibitem{\Power}
  {S. C. Power}
  {Simplicity of $C\sp\ast $-algebras of minimal dynamical systems}
  {{\it J. London Math. Soc. (2)}, {\bf 18} (1978), 534--538}

\bibitem{\Renault}
  {J. Renault}
  {A groupoid approach to C*-algebras}
  {Lecture Notes in Mathematics vol.~793, Springer, 1980}

\bibitem{\NewRenault}
  {J. Renault}
  {Cuntz-like algebras}
  {In {\it Operator theoretical methods}. Proceedings of the 17th
international conference on operator theory, Timisoara, Romania, June
23-26, 1998. Bucharest: The Theta Foundation, 2000, pp.~371--386}

\bibitem{\Stacey}
  {P. J. Stacey}
  {Crossed products of C*-algebras by *-endomorphisms}
  {{\it J. Aust. Math. Soc., Ser. A}, {\bf 54} (1993), 204--212}

\bibitem{\Takesaki}
  {M. Takesaki}
  {Theory of Operator Algebras I}
  {Springer-Verlag, 1979}

\bibitem{\Tom}
  {J. Tomiyama}
  {The Interplay between Topological Dynamics and Theory of
C*-algebras}
  {Lecture Notes Series vol.~2, Department of Mathematics, Seoul
National University, 1992}

\bibitem{\Watatani}
  {Y. Watatani}
  {Index for C*-subalgebras}
  {{\it Mem. Am. Math. Soc.}, {\bf 424} (1990), 117 p}
  \bye